\documentclass[12pt]{article}
\usepackage{fullpage,graphicx,psfrag,amsmath,amsfonts,verbatim,placeins}
\usepackage{algorithm,algpseudocode}
\usepackage[small,bf]{caption}
\usepackage{amssymb}
\graphicspath{{figures/}}

\usepackage{listings}
\usepackage{longtable}
\usepackage{subcaption}

\usepackage{color}

\definecolor{codegreen}{rgb}{0,0.6,0}
\definecolor{codegray}{rgb}{0.5,0.5,0.5}
\definecolor{codepurple}{rgb}{0.58,0,0.82}
\definecolor{backcolour}{rgb}{0.95,0.95,0.98}

\lstdefinestyle{mystyle}{
backgroundcolor=\color{backcolour},
keywordstyle=\color{magenta},
numberstyle=\tiny\color{codegray},
stringstyle=\color{codepurple},
basicstyle=\ttfamily\footnotesize,
breakatwhitespace=false,
breaklines=true,
captionpos=t,
keepspaces=true,
numbers=left,
numbersep=5pt,
showspaces=false,
showstringspaces=false,
showtabs=false,
tabsize=2,
}

\definecolor{seagreen}{rgb}{0.18, 0.55, 0.34}
\definecolor{mediumviolet-red}{rgb}{0.78, 0.08, 0.52}
\definecolor{violet}{rgb}{0.4, 0.0, 0.8}
\definecolor{khaki}{rgb}{0.94, 0.9, 0.55}

\lstdefinelanguage{mypython}
{
keywords=[1]{from, import, as, assert, not, print, nonneg, PSD, axis, lambda},
keywordstyle=[1]{\color{mediumviolet-red}},
keywords=[2]{cp, cpg, lo, pl, cvxpy, cvxpygen, Variable, Parameter, CallbackParam,
sqrt, exp, numpy, np, Problem, Minimize, Maximize, value, dual_value, solve, inner,
sum, multiply, arange, range, norm1, norm2, norm_inf, abs, square,
diagonal, outer, pos, hstack, power, generate_code},
keywordstyle=[2]{\color{violet}},
upquote=true,
showstringspaces=false,
basicstyle=\ttfamily\small,
columns=fullflexible,
keepspaces=true,
emph={True,False,def,return,float,class,match,switch,len, for, while, if, in, range, break},
emphstyle={\color{mediumviolet-red}},
belowskip=1em,
aboveskip=1em,
morecomment=[l]{\#}
}

\newcommand{\BEAS}{\begin{eqnarray*}}
\newcommand{\EEAS}{\end{eqnarray*}}
\newcommand{\BEA}{\begin{eqnarray}}
\newcommand{\EEA}{\end{eqnarray}}
\newcommand{\BEQ}{\begin{equation}}
\newcommand{\EEQ}{\end{equation}}
\newcommand{\BIT}{\begin{itemize}}
\newcommand{\EIT}{\end{itemize}}
\newcommand{\BNUM}{\begin{enumerate}}
\newcommand{\ENUM}{\end{enumerate}}

\newcommand{\BA}{\begin{array}}
\newcommand{\EA}{\end{array}}

\newcommand{\eg}{{\it e.g.}}
\newcommand{\ie}{{\it i.e.}}

\newcommand{\ones}{\mathbf 1}

\newcommand{\reals}{{\mbox{\bf R}}}
\newcommand{\integers}{{\mbox{\bf Z}}}

\newcommand{\diag}{\mathop{\bf diag}}

\newcommand{\Expect}{\mathop{\bf E{}}}
\newcommand{\Prob}{\mathop{\bf Prob}}

{\begin{quote}}{\end{quote}}

\makeatletter
\long\def\@makecaption#1#2{
\vskip 9pt
\begin{small}
\setbox\@tempboxa\hbox{{\bf #1:} #2}
\ifdim \wd\@tempboxa > 5.5in
\begin{center}
\begin{minipage}[t]{5.5in}
\addtolength{\baselineskip}{-0.95pt}
{\bf #1:} #2 \par
\addtolength{\baselineskip}{0.95pt}
\end{minipage}
\end{center}
\else
\hbox to\hsize{\hfil\box\@tempboxa\hfil}
\fi
\end{small}\par
}
\makeatother

\newcounter{oursection}

\newcounter{lecture}

\graphicspath{{img/}}

\usepackage{mathtools}
\usepackage{adjustbox}
\usepackage{booktabs}
\usepackage{csquotes}

\usepackage{hyperref}
\hypersetup{ colorlinks   = true,    
urlcolor     = blue,    
linkcolor    = black,    
citecolor    = blue      
}

\usepackage{color}

\title{Estimating Price Elasticity Matrices}
\author{Maximilian Schaller 
\and Stephen Boyd}
\date{April 12, 2026}

\begin{document}
\maketitle

\begin{abstract}
The relationship between demand and prices of a set of products can be modeled as a linear
mapping from logarithmic price changes to logarithmic changes in demand.
We consider the problem of estimating the coefficient matrix of this mapping, the elasticity
matrix, based on observed data consisting of real-valued prices and integer-valued demands.
We regularize the estimation problem by imposing a factor model structure,
\ie, that the elasticity matrix is diagonal plus low-rank, similar to factor models used
for financial returns.
Maximizing the likelihood of observations of this model is a bi-convex problem, 
meaning that there is a partition of the variables in which it is convex in 
each set when the other is fixed.
We propose and compare three methods for finding a locally optimal estimate. 
The first is based on alternating maximization, and
involves solving a sequence of convex problems. The second method exploits efficient
gradient computations in a gradient ascent method. The final method is to use a general purpose
nonlinear programming method. While all methods give the same result on numerical
examples, the gradient ascent
method is substantially faster, due to its efficient gradient evaluations.
We report the likelihood with different hyper-parameters for
synthetic and real-world data, with similar results. For synthetic data, we also report the
realized profit when using the elasticity estimate for optimal pricing, which is maximized for
the same set of hyper-parameters that also maximizes the likelihood.
This paper is accompanied by easy to use open source Python code for fitting 
elasticity matrices to observed data, using our three numerical methods.
\end{abstract}

\clearpage
\tableofcontents
\clearpage

\section{Introduction}

\subsection{Price elasticity of the demand}
Already in the 1960s, the relationship between demand
and prices of multiple products has been analyzed,
in terms of the Rotterdam model \cite{barten1964consumer,theil1965},
which is a local first-order approximation of utility-maximizing demand systems
\cite[\S3.D]{armstrong2018multiproduct,mas1995microeconomic}.
In particular, the Rotterdam model relates local changes in log-demand
linearly to local changes in log-prices. The coefficient matrix is referred
to as the \emph{elasticity matrix} or \emph{Slutsky substitution matrix}
\cite[\S2.F]{mas1995microeconomic}, and is constant. The diagonal entries
of such a matrix are typically referred to as \emph{self-elasticities} and
the off-diagonal entries are typically referred to as \emph{cross-elasticities}.
The Translog model \cite{christensen1975transcendental} from the 1970s
and the PIGL and PIGLOG architectures \cite{muellbauer1975aggregation},
including the widely used almost ideal demand system (AIDS) \cite{deaton1980almost}
and LA-AIDS \cite{deaton1980economics} from the 1980s, compute budget shares
as linear functions of log-price changes. One can convert those coefficients
to demand elasticities, which depend on the budget shares, and therefore on
the prices themselves.
The exact affine Stone index (EASI) \cite{lewbel2009tricks} model from the 2000s
is a generalization of AIDS, and shares the property of non-constant elasticity.

The above models generally result in full-rank elasticity matrices, potentially
rendering the estimation of such for large product portfolios computationally
expensive, or intractable.
The nested logit \cite{berry1994estimating} and nested constant elasticity of
substitution (CES) \cite{broda2006globalization} models from the 1990s and 2000s,
respectively, impose structure on the demand model,
by grouping products into nests. Within a nest, cross-elasticities are larger, and
across nests, cross-elasticities are smaller.
This structure corresponds to a low-rank plus diagonal elasticity matrix.
The low-rank part can be constructed by summing a rank-1 matrix per nest,
plus the own-price effects.
Recent work on low-rank cross-price effects \cite{fally2025consumer}
imposes similar structure, by using \emph{aggregators}, which are functions
that map information from the set of prices to a scalar, which
is then mapped to the demand for all products.

While these approaches share the idea of low-rank cross-elasticity,
motivated from first principles, we propose a generic
low-rank plus diagonal structure for identifying elasticity matrices
in a data-driven way.

\subsection{Our contribution}\label{s-contribution}
We provide a simple and efficient method for fitting low-rank plus diagonal
elasticity matrices to data, that can be used for, \eg, optimal pricing
\cite{wilson1993nonlinear,gallego1994optimal,ferreira2016analytics,schaller2025note}.

A low-rank plus diagonal elasticity matrix $E \in \reals^{n \times n}$ (for $n$ goods)
is one of the form
\[
E = BC^T + \diag(s),
\]
where the factors $B, C \in \reals^{n \times r}$ have rank $r < n$, and $s \in \reals^n$.
As described above, there are several interpretations of this structure
\cite{broda2006globalization,fally2025consumer}. We suggest a general,
statistical interpretation. Each column of $C$ is a principal direction
in the space of log-price changes. The columns of $B$ attribute
price movements along these directions, respectively,
to changes in the principle directions of log-demand. The part of the
self-elasticities that is not captured by the low-rank part is captured in $s$.

The same structure is widely used to model the (symmetric, positive definite) 
covariance of financial returns
\cite{johansson2023simple,lettau2020estimating,lettau2020factors,fama1992cross,fama1993common}.
In that case, $B$ and $C$ (which are equal in this case) are the \emph{factor loadings}, \ie,
mappings between $n$ financial returns and $r$ statistical factors
\cite{lettau2020estimating,lettau2020factors} or fundamental  factors
\cite{fama1992cross,fama1993common}
(that explain the returns), times a Cholesky root of the factor covariance.
The diagonal part $s$ (which is positive in this case) 
then corresponds to the idiosyncratic risk of the the individual financial instruments.

We explore several methods for fitting such elasticity matrices to data, where
the demand is measured in counts of goods sold, \ie, is integer-valued,
and prices are real-valued.  For typical measurement intervals of, say,
large online retailers, counts can be large (\eg, thousands or more),
so that one could reasonably approximate demand as real-valued, but we do not
make this assumption.

One method we use is \emph{alternating maximization}
\cite{jain2013low,udell2016generalized},
where in each iteration one of the factors $B$ and $C$ of the low-rank part of the matrix
is held constant, which yields a convex optimization problem that can be effectively solved.
The other method is a gradient ascent method with a simple step size schedule,
where we make heavy use of fast gradient computations that exploit the specific problem
structure.
We compare the two methods to a generic nonlinear programming solver.
We find that all methods give the same result on our numerical examples, so none produces
better fits.
Our gradient ascent method is, however, much faster than the other two methods,
due to its efficient gradient computations.

We evaluate the elasticity estimates statistically in terms of likelihood under a 
simple Poisson model.
In synthetic examples we also evaluate the results in terms of realized profit when 
used within the product pricing problem from \cite{schaller2025note}, 
which is a natural end use of an elasticity model.
It is re-assuring that the model hyper-parameters
(rank and the regularization factor for the size of the low-rank part)
are the same when chosen using these two evaluation methods.
For real data from the Dominick's Finer Foods (DFF) dataset \cite{dominicks1997},
the likelihood behaves similarly as a function of the hyper-parameters.

Code and data to reproduce the results of this paper, as well as to fit general 
elasticity matrices with different structure, is available at
\begin{center}
\url{https://github.com/cvxgrp/elasticity-estimation}.
\end{center}

\subsection{Outline}
In \S\ref{s-elasticity-matrices} we introduce a simple Poisson model
for the demand given prices of multiple products, and describe how it relates
to a log-linear demand model, where a price elasticity matrix
maps changes in log-prices to changes in log-demand.
In \S\ref{s-estimation} we describe how the maximum likelihood estimate
for such an elasticity matrix can be obtained via convex optimization,
and how the estimation problem becomes nonconvex when a
low-rank plus diagonal structure is imposed on the elasticity matrix.
In \S\ref{s-fitting} we present three numerical methods for fitting a
low-rank plus diagonal elasticity matrix to data.
We define evaluation metrics in \S\ref{s-evaluation} before
reporting experimental results in \S\ref{s-results}, for synthetic
and real demand data.

\section{Price elasticity matrices}\label{s-elasticity-matrices}
We follow the notation from \cite{schaller2025note}.
For a collection of $n$ goods (products, services, etc.), the demand $d_i$ is the
(integer) quantity of product $i$ sold at price $p_i$ in a given 
time interval.
The respective vectors of demand and price are $d\in \integers^n_+$
and $p \in \reals_+^n$.

\subsection{Poisson model of the demand}\label{s-poisson}
We model the demand, which takes on integer values, as an IID 
Poisson distribution,
\cite{feller1991introduction} (which is often used in demand modeling
\cite{silva2006log,burda2012poisson}),
\BEQ\label{e-poisson-dist}
\Prob(d_i = k) = \frac{\lambda_i(p)^k e^{-\lambda_i(p)}}{k!},
\quad k = 0,1,\ldots, \quad i=1,\ldots,n,
\EEQ
where the rates $\lambda_i \in \reals_{++}$ 
(equal to $\Expect d_i$) are functions of the price vector $p$.
In particular, we use the Poisson regression
model \cite{coxe2009analysis,hayat2014understanding} with $n$ features
\BEQ\label{e-poisson-model}
\log \lambda_i = \alpha_i + \sum_{j=1}^n \beta_{ij} f_j(p_j), \quad i=1,\ldots,n,
\EEQ
where $\alpha_i \in \reals$ and $\beta_{ij} \in \reals$ are model parameters
and $f_j : \reals_+ \rightarrow \reals$ is a function of the $j$th price.
This is sometimes referred to as a generalized linear model~\cite{hastie2017generalized}
with a logarithmic link function.
Here, the $j$th feature $f_j(p_j)$ depends on the $j$th price $p_j$.

\subsection{Price elasticity of the demand}
We consider (positive) nominal prices $p^\text{nom}_1, \ldots, p^\text{nom}_n$,
which are typically the prices at which the products
(or similar reference products, respectively) have been sold in the past.
We assume that they are known.
We denote the demand for product $i$ at its nominal price $p_i^\text{nom}$
as $d_i^\text{nom}$, which is not known, in general, and real-valued as
it represents the respective Poisson rate.

In the model~\eqref{e-poisson-model}, we may interpret
the offset $\alpha_i$ as $\log d_i^\text{nom}$
and set
\[
f_j(p_j) = \log(p_j / p^\text{nom}_j),
\]
so we can re-write~\eqref{e-poisson-model} as
\[
\log(\Expect d_i / d^\text{nom}_i) = \sum_{j=1}^n \beta_{ij} \log(p_j / p^\text{nom}_j).
\]
We see that it is natural to view the relationship between demand and prices
in terms of the logarithmic fractional (expected) demand and price changes
\[
\delta_i = \log (\Expect d_i / d_i^\text{nom}), \quad
\pi_i = \log (p_i / p_i^\text{nom}), \quad i=1,\ldots,n,
\]
respectively. We let $\delta, \pi \in \reals^n$ denote the vectors of demand and price changes,
respectively.
Using matrix notation, we arrive at
\[
\delta = E \pi,
\]
where $E \in \reals^{n \times n}$ is the elasticity matrix as introduced in \S\ref{s-contribution},
with $E_{ij} = \beta_{ij}$, the elasticity of the demand for product $i$
with respect to a changed price of product $j$.

\section{Estimating elasticity matrices}\label{s-estimation}

\subsection{Maximum likelihood estimate}\label{s-mle}
We consider a sample of $N$ observations $d^{(j)} \in \integers^n_+, p^{(j)} \in \reals^n_+$
drawn IID from \eqref{e-poisson-dist}.
We denote the respective
price changes from nominal by $\pi^{(j)} = \log(p^{(j)} / p^\text{nom})$,
and the rates according to~\eqref{e-poisson-model} by $\lambda^{(j)}$.
The log-likelihood takes the form
\[
\sum_{j=1}^N \sum_{i=1}^n ( d^{(j)}_i \log \lambda^{(j)}_i - \lambda^{(j)}_i ),
\]
where we dropped constant terms, and the dependence on the model parameters
is implicit via the Poisson regression model~\eqref{e-poisson-model}.

We define $y^{(j)}_i = \log \lambda^{(j)}_i$ and
\[
\tilde E = \left[\begin{array}{cc} E & \log d^\text{nom} \end{array}\right], \quad
\tilde \pi^{(j)} = (\pi^{(j)}, 1), \quad j=1,\ldots,N,
\]
where $\log$ is taken elementwise, so $y^{(j)} = \tilde E \tilde \pi^{(j)}$.
We arrange the observed sales and price changes as
\[
D = \left[\begin{array}{ccc}
d^{(1)} & \cdots & d^{(N)}
\end{array}\right], \quad
\tilde \Pi = \left[\begin{array}{ccc}
\tilde \pi^{(1)} & \cdots & \tilde \pi^{(N)}
\end{array}\right].
\]
To obtain the maximum likelihood estimate (MLE) for $E$ and $d^\text{nom}$,
we solve
\BEQ\label{e-mle}
\begin{array}{ll}
\mbox{maximize} & f(\tilde E) =
\frac{1}{N} \sum_{j=1}^N \sum_{i=1}^n (D \circ \tilde E \tilde \Pi - \exp(\tilde E \tilde \Pi))_{ij}
\end{array}
\EEQ
where $\tilde E \in \reals^{n \times n+1}$ is the variable and exp is taken elementwise,
and $\circ$ denotes the Hadamard or elementwise product.
We extract the MLE for $E$ and $d^\text{nom}$ from the solution $\tilde E^\star$.
This is a convex optimization problem, as the first part of the objective as affine,
and the second is concave \cite{boyd2004convex}.

\subsection{Low-rank plus diagonal elasticity}
In typical consumer data, the squared number of goods exceeds the number of data points,
so estimating a full-rank $n$ by $n$ elasticity matrix is not tractable.
Instead, we consider the low-rank plus diagonal structure introduced in \S\ref{s-contribution},
\BEQ\label{e-low-rank-plus-diag}
E = BC^T + \diag(s),
\EEQ
where $B, C \in \reals^{n \times r}$ and $s \in \reals^n$, so the number of
estimated parameters is $2nr + n = n (2r + 1)$. When $r \ll n$, this is much less than the $n^2$
parameters of a full-rank elasticity matrix.
Imposing this structure on the elasticity matrix makes parameter estimation tractable,
and more importantly improves the out-of-sample prediction of demand from prices.

\paragraph{Low-rank plus diagonal MLE.}
We add~\eqref{e-low-rank-plus-diag} as a constraint
and a quadratic regularizer to the MLE problem~\eqref{e-mle}, and solve
\BEQ\label{e-mle-low-rank}
\begin{array}{ll}
\mbox{maximize} &
\frac{1}{N} \sum_{j=1}^N \sum_{i=1}^n (D \circ \tilde E \tilde \Pi - \exp(\tilde E \tilde \Pi))_{ij}
- \frac{\lambda}{2} ( \|B\|_F^2 + \|C\|_F^2 ) \\
\mbox{subject to} & \tilde E_{1:n} = BC^T + \diag(s).
\end{array}
\EEQ
Here, the variables are $\tilde E \in \reals^{n \times n + 1}$
(the subscript $1:n$ denotes the first $n$ columns), $B, C \in \reals^{n \times r}$,
and $s \in \reals^n$. The value of $\lambda > 0$ controls the regularization strength.
In fact, this regularizer has an equivalent effect to the dual spectral norm penalty
$\lambda \|M\|_*$ (the sum of singular values),
when imposing $E = M + \diag(s)$ with $\textbf{rank}(M) \leq r$
\cite{udell2016generalized}. Therefore, we can view $r$ as an upper bound on
the rank of the resulting low-rank part, and $\lambda$ as a tuning knob for the
resulting rank.

This is not a convex optimization problem, since the equality constraint involves the
nonlinear expression $BC^T$. However, when either $B$ or $C$ are fixed, the problem
becomes convex \cite{boyd2004convex}, which we will exploit.

\section{Fitting methods}\label{s-fitting}

We describe three methods for solving problem \eqref{e-mle-low-rank}.
In \S\ref{s-am}, we describe alternating maximization, where we alternate between
fixing $B$ or $C$, and solving a respective  problem that is convex
in the remaining variables.
In \S\ref{s-ga}, we derive an expression for the gradient of the objective of
problem \eqref{e-mle-low-rank} with respect to $B$, $C$, $s$, and $\log d^\text{nom}$,
which can be evaluated efficiently, in a vectorized way.
The third method is described in \S\ref{s-nlp} and applies a general nonlinear programming
solver to problem \eqref{e-mle-low-rank} as is, with variables $\tilde E$, $B$, $C$, and $s$.

\subsection{Alternating maximization}\label{s-am}
Problem~\eqref{e-low-rank-plus-diag} becomes convex when fixing either
$B$ or $C$, since then the constraint becomes an affine equality constraint (and the
maximized objective was concave already) \cite{boyd2004convex}.
Therefore, we consider an alternating maximization scheme
where we alternate between computing solutions to variants of \eqref{e-mle-low-rank}
with $B$ or $C$ fixed, respectively.

We initialize $\hat C$ randomly, away from zero.
Then, we take the variant of \eqref{e-mle-low-rank} with $C$ declared a parameter,
assign $\hat C$ to $C$ and obtain a solution with a convex optimization solver,
from which we take $\hat B$.
Next, we take the variant of \eqref{e-mle-low-rank} with $B$ declared a parameter,
assign $\hat B$ to $B$, solve the problem, and obtain $\hat C$ from the solution.
We repeat this procedure until the objective value converges.

This alternating maximization scheme is an ascent method.
Since each sub-problem is a convex maximization problem, the respective optimal
value can be computed effectively, and will be at least as high as the value at the
initialization (from the previous iteration).

\paragraph{CVXPY specification.}
Figure~\ref{f-cvxpy-am} shows
how the procedure is implemented in the domain-specific modeling
language for convex optimization CVXPY~\cite{diamond2016cvxpy,agrawal2018rewriting}.
In lines 4--22, the two respective problems are modeled with CVXPY.
In lines 10 and 11, \verb|Parameter| objects are declared to
allow for efficient updates to the respective problems when the values of $B$
and $C$ change. 
In lines 25--30, five iterations of alternating maximization are run.

\begin{figure}
\lstset{language=mypython,
numbers=left,
xleftmargin=1.5em,
linewidth=0.96\columnwidth}
\begin{lstlisting}[frame=lines]
import cvxpy as cp
	
# variables
Etilde = cp.Variable((n, n + 1))
B = cp.Variable((n, r))
C = cp.Variable((n, r))
s = cp.Variable(n)

# parameters
Bp = cp.Parameter((n, r))
Cp = cp.Parameter((n, r))

# objectives and constraints
f = cp.sum(cp.multiply(D, Etilde @ Pitilde) - cp.exp(Etilde @ Pitilde)) / N
obj_over_B = cp.Maximize(f - (lam / 2) * cp.sum_squares(B))
obj_over_C = cp.Maximize(f - (lam / 2) * cp.sum_squares(C))
constr_over_B = [Etilde[:, :-1] == B @ Cp.T + cp.diag(s)]
constr_over_C = [Etilde[:, :-1] == Bp @ C.T + cp.diag(s)]

# problems
prob_over_B = cp.Problem(obj_over_B, constr_over_B)
prob_over_C = cp.Problem(obj_over_C, constr_over_C)

# initialize and solve
Cp.value = ...
for _ in range(5):
    prob_over_B.solve()
    Bp.value = B.value
    prob_over_C.solve()
    Cp.value = C.value
\end{lstlisting}
\caption{Alternating maximization with CVXPY.
The dimensions \texttt{n}, \texttt{r}, \texttt{N},
and the data \texttt{D}, \texttt{Pitilde}, \texttt{lam} are given.}
\label{f-cvxpy-am}
\end{figure}

\subsection{Gradient ascent}\label{s-ga}
We run a simple gradient ascent method,
where we take steps along the gradients of the objective in~\eqref{e-mle-low-rank},
denoted by $g(X)$, with respect to
$X = \left[\begin{array}{cccc}
B & C & s & \log d^\text{nom}
\end{array}\right] \in \reals^{n \times 2r + 2}$.
We use a simple line search to guarantee that the algorithm is an ascent method.
If the likelihood is improved with a gradient step of the current step size
(denoted by $\alpha$),
we use it for the current iteration
and increase it by a constant factor $\gamma > 1$ for the next iteration.
Otherwise, we repeatedly shrink the step size by a constant factor $\eta > 1$
until the likelihood is improved.
The method is given in algorithm~\ref{alg-ga}.

\begin{algorithm}
\caption{Gradient ascent}\label{alg-ga}
\begin{algorithmic}[1]
\State Initialize $X^0$, $\alpha^0$, $k=0$
\Repeat
\State $\hat X = X^k + \alpha^k \nabla g(X)$
\Comment tentative update
\If{$g(\hat X) > g(X^k)$}
\State $X^{k+1} = \hat X$, $\alpha^{k+1} = \gamma \alpha^k$
\Comment accept update and increase step size
\Else
\State $\alpha^k \gets \alpha^k / \eta$, go to step 3
\Comment shrink step size and re-evaluate
\EndIf
\State $k \leftarrow k + 1$
\Until{$\| \nabla g(X^k)\|_F \leq \epsilon^\text{rel} \|X^k\|_F + \epsilon^\text{abs}$}
\end{algorithmic}
\end{algorithm}

\paragraph{Initialization.}
We draw the entries of $B^0$ and $C^0$ IID from $\mathcal{N}(0, 1 / (n\sqrt{r}))$,
set $s^0 = -\ones$, and $(\log d^\text{nom})^0 = 0$.
The algorithm is not particularly dependent on the initial step size $\alpha^0$
and the line search parameters $\gamma$ and $\eta$.
Reasonable choices are, \eg, $\alpha^0 = 1.0$, $\gamma = 1.2$, and $\eta = 1.5$.

\paragraph{Gradient computation.}
The gradient of the objective in \eqref{e-mle-low-rank} with respect to $\tilde E$ is
\[
\Delta = (1/N) (D - \exp(\tilde E \tilde \Pi)) \tilde \Pi^T
= (1/N) (D  \tilde \Pi^T -  \exp(\tilde E \tilde \Pi) \tilde \Pi^T).
\]
The second expression is faster to compute for changing values of $\tilde E$
when $N > n$, since $D \tilde \Pi^T$ needs to be computed only once ahead of time.
We apply the chain rule to obtain the gradient with respect to $X$,
\[
\nabla g(X) = \left[\begin{array}{cccc}
\Delta_{1:n} C - \lambda B &
\Delta_{1:n}^T B - \lambda C &
\diag(\Delta_{1:n}) &
\Delta_{n+1}
\end{array}\right],
\]
where subscript $n+1$ denotes column $n+1$.
We can compute $\nabla g(X)$ very fast, as it involves mostly
BLAS level 3 computations \cite{dongarra1990set}.

\paragraph{Stopping criterion.}
We stop the algorithm as soon as the termination criterion
\[
\| \nabla g(X^k)\|_F \leq \epsilon^\text{rel} \|X^k\|_F + \epsilon^\text{abs}
\]
with $\epsilon^\text{rel}, \epsilon^\text{abs} > 0$ is met.
Typical values for $\epsilon^\text{rel}$
and $\epsilon^\text{abs}$ range between $10^{-6}$ and $10^{-2}$.

\paragraph{Python specification.}
Figure~\ref{f-python-ga} shows a Python implementation of
algorithm~\ref{alg-ga}. The variable, gradient, and step size
are initialized in lines 3--6, and the actual algorithm is
implemented via two simple loops in lines 8--17.

\begin{figure}
\lstset{language=mypython,
numbers=left,
xleftmargin=2em,
linewidth=0.96\columnwidth}
\begin{lstlisting}[frame=lines]
import numpy.linalg as la
	
X = ...
grad_X = grad(X)
g_values = [g(X)]
alpha = 1.0

while la.norm(grad_X) > eps_rel * la.norm(X) + eps_abs:
    while True:
        X_hat = X + alpha * grad_X
        if g(X_hat) > g_values[-1]:
            X = X_hat
            grad_X = grad(X)
            g_values.append(g(X))
            alpha *= 1.2
            break
        alpha /= 1.5
\end{lstlisting}
\caption{Gradient ascent in Python.
The initial value of \texttt{X}, the Python functions \texttt{g} and \texttt{grad},
and the tolerances \texttt{eps\_rel} and \texttt{eps\_abs} are given.}
\label{f-python-ga}
\end{figure}

\subsection{Nonlinear programming}\label{s-nlp}
One can view the constraint in problem~\eqref{e-mle-low-rank} as an instance of
general nonlinear (and differentiable) constraints and apply nonlinear programming
(NLP) methods, which are local methods as they use gradients or (approximate) Hessians
\cite{bertsekas1997nonlinear, liu1989limited, nocedal2006numerical, kuhn2013nonlinear}.
Well-known NLP solver implementations are the open-source IPOPT \cite{wachter2006implementation}
and the proprietary KNITRO \cite{byrd2006knitro}.

\paragraph{CVXPY specification.}
Figure~\ref{f-cvxpy-nlp} shows how problem~\eqref{e-mle-low-rank} is solved
via the CVXPY NLP interface~\cite{cederberg2025disciplined}.
The nonconvex problem is modeled in lines 4--13 with CVXPY.
After variable initialization in line 15, the problem is solved in line 16,
where \verb|nlp=True| informs CVXPY that the problem should be treated as
an NLP.

\begin{figure}
\lstset{language=mypython,
numbers=left,
xleftmargin=2em,
linewidth=0.97\columnwidth}
\begin{lstlisting}[frame=lines]
import cvxpy as cp

# variables
Etilde = cp.Variable((n, n + 1))
B = cp.Variable((n, r))
C = cp.Variable((n, r))
s = cp.Variable(n)

# objective, constraint, and problem
f = cp.sum(cp.multiply(D, Etilde @ Pitilde) - cp.exp(Etilde @ Pitilde)) / N
obj = cp.Maximize(f - (lam / 2) * (cp.sum_squares(B) + cp.sum_squares(C)))
constr = [Etilde[:, :-1] == B @ C.T + cp.diag(s)]
prob = cp.Problem(obj, constr)

# solve problem
Etilde.value, B.value, C.value, s.value = ...
prob.solve(nlp=True)
\end{lstlisting}
\caption{Nonlinear programming with CVXPY.
The dimensions \texttt{n}, \texttt{r}, \texttt{N},
and the data \texttt{D}, \texttt{Pitilde}, \texttt{lam} are given.}
\label{f-cvxpy-nlp}
\end{figure}

\section{Evaluation}\label{s-evaluation}
There are several ways to assess the quality of an elasticity estimate.
We focus on two different types of metrics: the first statistical, based on 
out-of-sample log-likelihood, and the second based on a typical end-use, realized profit
after optimally pricing the products using the elasticity model.
In the first, we compute the log-likehood of the observed demand under the
Poisson model with the estimated price elasticity matrix (and nominal demand).
This provides a direct measure of how well the estimator captures the
statistical structure of the demand.
In the second, we solve a profit maximization problem given the elasticity estimate,
and record the realized profit, given a simulated true demand function.

These metrics offer complementary insights. 
Log-likelihood quantifies how well the estimator aligns with the observed
pattern of price and demand changes.
However, it relies on the assumption that the demand follows
a Poisson distribution, which may or may not hold in practice.
Assessing the elasticity estimate in terms of profit does not require such
a distributional assumption and instead evaluates the economic usefulness
of the estimator by measuring its profit impact on pricing decisions.
At the same time, the profit after pricing also depends on
design choices, such as limits on price changes, making attribution less direct,
not to mention that we compute it only in simulation.
Fortunately we will see that these two different evaluations of proposed models
lead to very similar conclusions.

\subsection{Log-likelihood}
The average log-likelihood is
\BEQ\label{e-log-likelihood}
\mathcal{L}(\tilde E; D, \tilde \Pi) = 
\frac{1}{N} \sum_{j=1}^N \sum_{i=1}^n
(D \circ \tilde E \tilde \Pi - \exp(\tilde E \tilde \Pi) - \log (D !))_{ij}
\EEQ
where log and factorial are taken elementwise and we consider data
where $D$ contains
positive integers (counts of goods sold).

We perform $K$-fold cross-validation. We split the data $D, \tilde \Pi$ into $K$
equally sized column slices, or folds,
\[
D = \left[\begin{array}{ccc}
D^{(1)} & \cdots & D^{(K)}
\end{array}\right], \quad
\tilde\Pi = \left[\begin{array}{ccc}
\tilde\Pi^{(1)} & \cdots & \tilde\Pi^{(K)}
\end{array}\right].
\]
We then fit $K$ models $\tilde E^{(1)}, \ldots, \tilde E^{(K)}$,
leaving out one fold each time, for which we then compute the
respective log-likelihood. We report the average
log-likelihood over all folds,
\BEQ\label{e-log-likelihood-cv}
\frac{1}{K} \sum_{k=1}^K \mathcal{L}(\tilde E^{(k)}; D^{(k)}, \tilde \Pi^{(k)}).
\EEQ

\subsection{Pricing performance}
For a given estimate $\tilde E$, we extract $E$ and $d^\text{nom}$, from
which we compute the nominal revenues
$r^\text{nom} = d^\text{nom} \circ p^\text{nom}$
and nominal costs $\kappa^\text{nom} = d^\text{nom} \circ c$, given the
nominal prices $p^\text{nom}$ and the vector of cost per unit $c$.
We then solve the (nonconvex) profit maximization problem
from \cite[\S2.4]{schaller2025note}
\BEQ\label{e-pricing}
\begin{array}{ll} \mbox{maximize} & 
\sum_{i=1}^n \left( r_i^\text{nom} \exp(\delta_i + \pi_i) - 
\kappa_i^\text{nom} \exp(\delta_i)\right)\\
\mbox{subject to} & \delta = E\pi, \quad \pi^\text{min} \preceq \pi \preceq \pi^\text{max},
\end{array}
\EEQ
where $\delta, \pi \in \reals^n$ are the variables and
$\pi^\text{min}, \pi^\text{max} \in \reals^n$
are given limits on the price changes.

Again, we perform $K$-fold cross-validation.
We solve problem~\eqref{e-pricing}
$K$ times, with the data taken from $\tilde E^{(1)}, \ldots, \tilde E^{(K)}$,
respectively.
From the respective solution $\pi^{(k)}$, we simulate the profit
\[
\mathcal{P}(\pi^{(k)}) =
\sum_{i=1}^n \left( r^\text{nom,syn} \circ \exp((E^\text{syn} + I) \pi^{(k)}) - 
\kappa^\text{nom,syn} \circ \exp(E^\text{syn} \pi^{(k)})\right)_i,
\]
where superscript syn denotes synthetic data, \ie, we take a synthetic
elasticity matrix $E^\text{syn}$ and nominal demand $d^\text{nom,syn}$
(see \S\ref{s-synthetic}),
and set $r^\text{nom,syn} = d^\text{nom,syn} \circ p^\text{nom}$,
$\kappa^\text{nom,syn} = d^\text{nom,syn} \circ c$.
We report the average over all folds,
\BEQ\label{e-pp}
\frac{1}{K} \sum_{k=1}^K \mathcal{P}(\pi^{(k)}).
\EEQ

\section{Results}\label{s-results}
We fit elasticity matrices by applying alternating maximization,
gradient ascent, and nonlinear programming, as described in \S\ref{s-fitting}.
In gradient ascent, we use $\gamma = 1.2$ and $\eta = 1.5$ for the line search
and termination tolerances $\epsilon^\text{rel} = \epsilon^\text{abs} = 10^{-3}$.
We use the same termination tolerances for alternating maximization
and nonlinear programming. The results are not particularly dependent on
the tolerance value.
To solve the convex subproblems of alternating maximization, we use
the convex optimization solver MOSEK \cite{mosek}.
For nonlinear programming, we use the open-source NLP solver IPOPT
\cite{wachter2006implementation}.
We interface with both solvers via CVXPY
and use their respective default settings.
We run the experiments on an Apple M1 Pro.

\subsection{Synthetic data}\label{s-synthetic}

\paragraph{Data set.}
We generate the entries of $p^{(1)}, \ldots, p^{(N)}$ uniformly
and IID from $[1, 2]$ and set
\[
\log p^\text{nom} = (1/N) \sum_{j=1}^N \log p^{(j)},
\]
where $\log$ is taken elementwise. Then, we
compute the log price changes
\[
\pi^{(j)} = \log p^{(j)} -\log p^\text{nom}, \quad j=1,\ldots,N
\]
and assemble them in $\tilde \Pi$ as described in \S\ref{s-mle}.

We construct a synthetic elasticity matrix
\[
E^\text{syn} = B^\text{syn} (C^\text{syn})^T + \diag(s^\text{syn}),
\]
where the entries of $B^\text{syn}, C^\text{syn} \in \reals^{n \times r^\text{syn}}$
are generated IID from $\mathcal{N}(0, 0.1)$, and the entries of $s^\text{syn}$
are generated uniformly and IID from $[-5, -1]$.

Ultimately, we set $d^\text{nom,syn} = \ones$ so
$\log d^\text{nom,syn} = 0$ and the Poisson rates are
\[
\lambda^{(j)} = \exp(E^\text{syn} \pi^{(j)}), \quad j=1,\ldots,N.
\]
We then generate the columns $d^{(j)}$ of $D$ from Poisson
distributions with respective rates $\lambda^{(j)}$.
We generate the entries of the cost vector $c$ IID from $[0.8, 1.2]$.

\paragraph{Results.}
For $n=100$, $N=200$, and $r^\text{syn} = 10$, we compare the
time it takes to solve problem \eqref{e-mle-low-rank} where $r = 10$ and
$\lambda = 0.1$, with alternating maximization (AM), gradient ascent,
and nonlinear programming (NLP), as described in \S\ref{s-fitting}.
All three methods converge to the same objective value of $-27.1$
(up to the relative tolerance).
Table \ref{t-times} shows the respective solve times.
\begin{table}[h]
\centering
\begin{tabular}{c|c|c}
\textbf{AM} & \textbf{Gradient ascent} & \textbf{NLP} \\
\hline
230 s & 0.202 s & 49.2 s \\
\end{tabular}
\caption{Solve times with AM, gradient ascent, and NLP.}
\label{t-times}
\end{table}
Gradient ascent is two to three orders of magnitude faster
than NLP and AM, respectively, due to its use of efficient (vectorized) gradient
computations. The other two mehods have much overhead, \eg, many
exponential cones when using AM \cite{mosek}.
We use gradient ascent for the remaining analyses.

\begin{figure}
\centering
\begin{subfigure}{0.48\columnwidth}
\centering
\includegraphics[width=\linewidth]{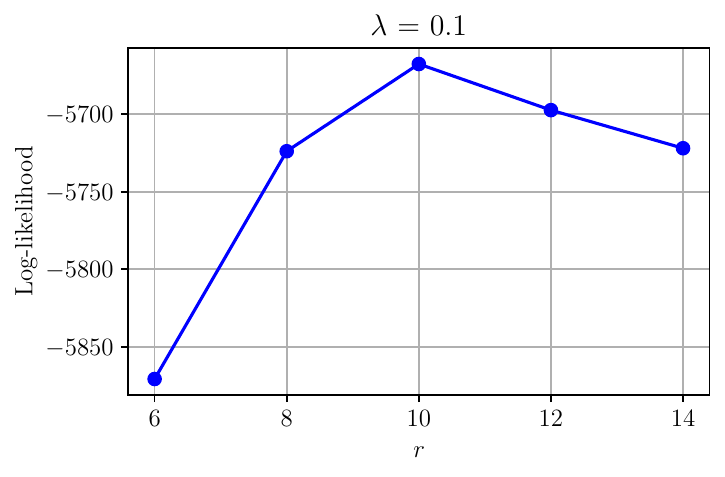}
\caption{Log-likelihood for varying $r$ at $\lambda=0.1$.}
\label{f-log-likelihood-syn-r}
\end{subfigure}
\hfill
\begin{subfigure}{0.48\columnwidth}
\centering
\includegraphics[width=\linewidth]{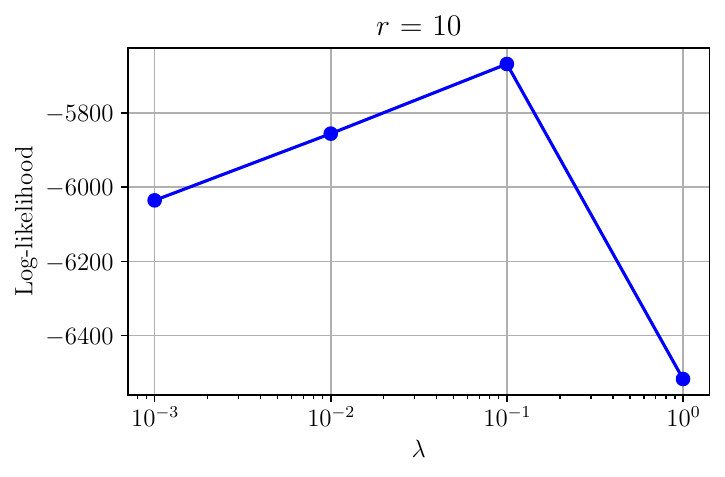}
\caption{Log-likelihood for varying $\lambda$ at $r=10$.}
\label{f-log-likelihood-syn-lambda}
\end{subfigure}
\caption{Log-likelihood for varying $r$ and $\lambda$, with synthetic data.}
\label{f-log-likelihood-syn}
\end{figure}

\begin{figure}
\centering
\includegraphics[width=0.95\linewidth]{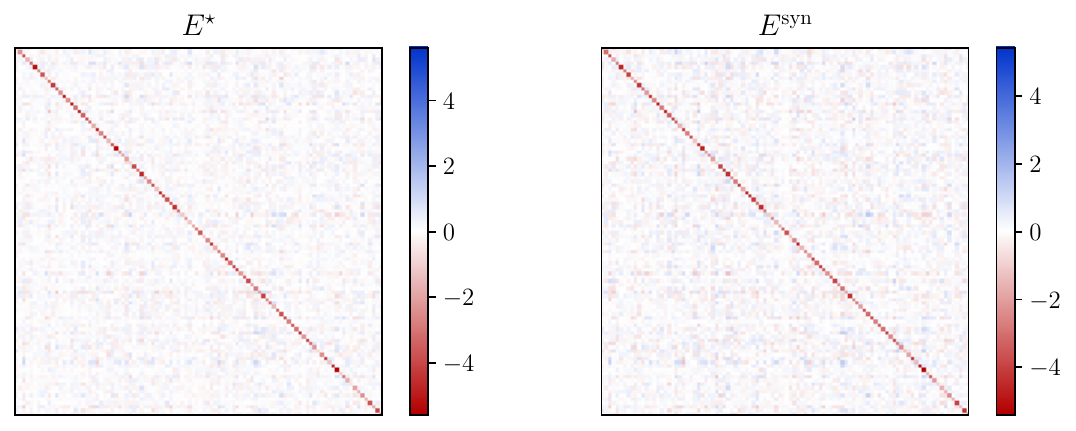}
\caption{Estimated elasticity matrix $E^\star$ for $\lambda=0.1$ and $r = 10$
and true (synthetic) elasticity matrix $E^\text{syn}$.}
\label{f-matrix-syn}
\end{figure}

We cross-validate
our elasticity estimator obtained by solving problem \eqref{e-mle-low-rank}
for
\[
r \in \{6, 8, 10, 12, 14\}, \quad \lambda \in \{ 10^{-3}, 10^{-2}, 10^{-1}, 10^0\}.
\]
We report the cross-validated log-likelihood~\eqref{e-log-likelihood-cv}
for $K=5$ folds. We obtain the highest log-likelihood at $r=10$
(equaling $r^\text{syn}$, as expected) and $\lambda=0.1$.
Figure~\ref{f-log-likelihood-syn} shows the log-likelihood over $r$ and $\lambda$,
respectively.
For the optimal configuration, figure~\ref{f-matrix-syn} compares
the elasticity estimate $E^\star$ to the true (synthetic) elasticity matrix
$E^\text{syn}$. We observe that our method recovers the true elasticity matrix,
up to slightly shrunk cross-elasticities due to regularization.

\begin{figure}
\centering
\begin{subfigure}{0.48\columnwidth}
\centering
\includegraphics[width=\linewidth]{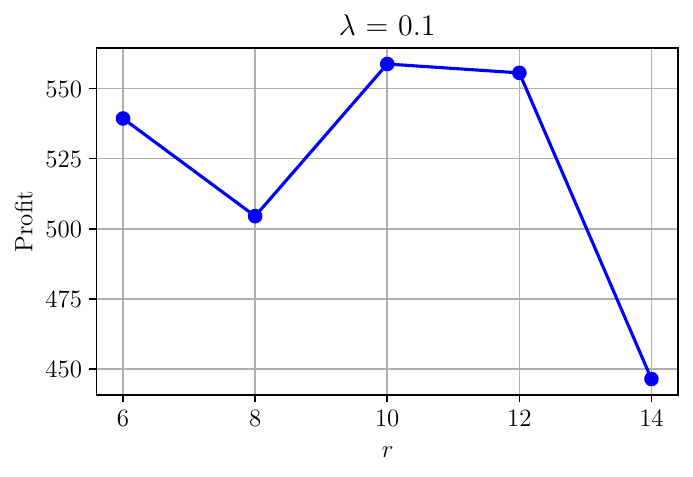}
\caption{Simulated profit for varying $r$ at $\lambda=0.1$.}
\label{f-pricing-performance-r}
\end{subfigure}
\hfill
\begin{subfigure}{0.48\columnwidth}
\centering
\includegraphics[width=\linewidth]{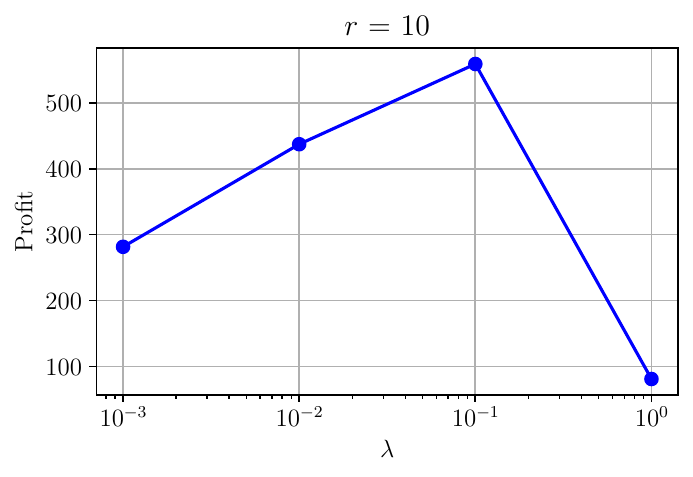}
\caption{Simulated profit for varying $\lambda$ at $r=10$.}
\label{f-pricing-performance-lambda}
\end{subfigure}
\caption{Pricing performance for varying $r$ and $\lambda$, with synthetic data.}
\label{f-pricing-performance}
\end{figure}

We also compute the cross-validated pricing performance~\eqref{e-pp}, \ie,
the simulated profit after solving the optimal pricing problem~\eqref{e-pricing}
with our elasticity estimate. We take $\pi^\text{min} = \log(0.8) \ones$ and
$\pi^\text{max} = \log(1.2) \ones$ in problem~\eqref{e-pricing}, \ie, we allow price
changes between $\pm 20\%$.
Consistent with the log-likelihood,
the highest pricing performance is obtained for $r=10$ and $\lambda=0.1$,
as shown in figure~\ref{f-pricing-performance}. The maximum cross-validated
profit is $558$, compared to a mean of $662$
when solving problem~\eqref{e-pricing}
with the true elasticity matrix (and true nominal demand).

\subsection{Real data}

\paragraph{Data set.}
We use the Dominick's Finer Foods (DFF) dataset provided by the
Kilts Center for Marketing at the University of Chicago Booth School of Business
\cite{dominicks1997}. The dataset contains scanner data from purchases
at DFF stores in the Chicago area from the 1990s, aggregated to weekly time intervals.
Unlike more recent datasets, the DFF dataset is publicly available to anyone.

We consider $n=20$ of the most frequently sold beverages, of which $N=360$
complete data points (out of $400$) are available;
we only consider weeks where all 20 beverages have been sold, and therefore have a known price.
(In scanner data, the price during a certain period is unknown when the product
is not sold in that period.)
We accept this selection bias as we demonstrate our numerical fitting method.
We take the nominal log-prices as the average of the log-prices.

\begin{figure}
\centering
\begin{subfigure}{0.48\columnwidth}
\centering
\includegraphics[width=\linewidth]{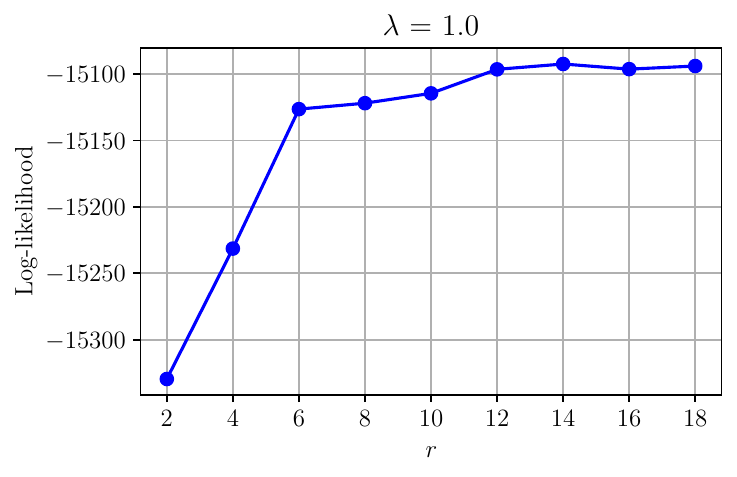}
\caption{Log-likelihood for varying $r$ at $\lambda=1.0$.}
\label{f-log-likelihood-real-r}
\end{subfigure}
\hfill
\begin{subfigure}{0.48\columnwidth}
\centering
\includegraphics[width=\linewidth]{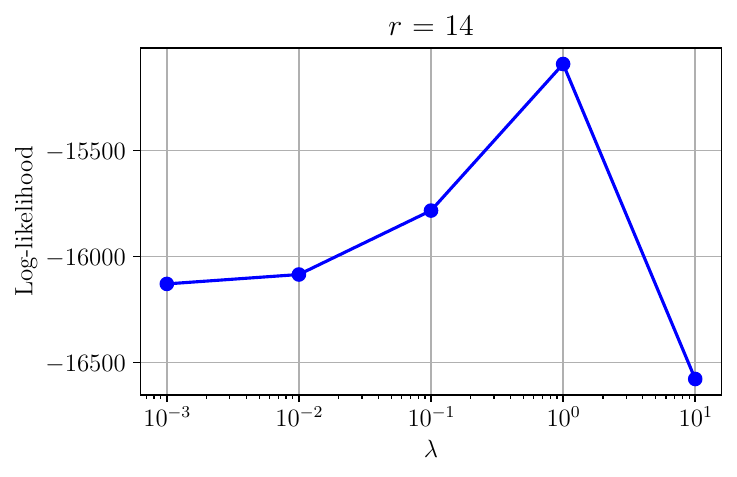}
\caption{Log-likelihood for varying $\lambda$ at $r=14$.}
\label{f-log-likelihood-real-lambda}
\end{subfigure}
\caption{Log-likelihood for varying $r$ and $\lambda$, with real data.}
\label{f-log-likelihood-real}
\end{figure}

\paragraph{Results.}
We consider $r \in \{2, 4, \ldots, 18\}$ and
$\lambda \in \{ 10^{-3}, 10^{-2}, \ldots, 10^1 \}$.
The cross-validated log-likelihood~\eqref{e-log-likelihood-cv} is shown
for $K=5$ folds in figure~\ref{f-log-likelihood-real}.
We obtain the highest log-likelihood at $r=14$ and $\lambda=1.0$.
The median relative error of the predicted (mean) demand with respect
to the true demand, over all products and time steps, is $34\%$
(cross-validated over the same five folds).
The median relative error due to the (predicted) Poisson noise is $12\%$,
\ie, about one third of the overall error.

\section*{Acknowledgments}
The authors thank Alexander Thebelt,  Marco Tacke, and Lutz Gruber for inspiring this work.

\clearpage
\bibliography{refs}

\end{document}